\documentstyle [12pt] {article}

\thicklines                         
\def\be{\begin{equation}}
\def\ee{\end{equation}}

\begin{document}

\par
\begin{center}
\vspace{0.3in}
\ \\ [0.3in]
{\Large {\bf $h$-deformation of  $GL(1|1)$}}\\
   [.9in]
Ludwik D\c abrowski and Preeti Parashar\\ [.2in]
{\it SISSA, Via Beirut 2-4, 34014 Trieste, Italy}\\ [.1in]
{ABSTRACT}\\[.3in]
\end{center}
\begin{quotation}
$h$-deformation of (graded) Hopf algebra of functions on supergroup $GL(1|1)$
is introduced via a contration of $GL_q (1|1)$.
The deformation parameter $h$ is odd (grassmann).
Related differential calculus on $h$-superplane is presented.\\
[.8in]
\end{quotation}
\vskip1cm
\centerline{Ref. SISSA: 137/95/FM}
\pagebreak

Recently a new class of quantum deformations of Lie groups and Lie algebras,
known as an `h-deformation', has been intensively studied [1-8, 17].
These deformations are important for at least two reasons.
One is that the deformation parameter $h$ is naturally dimensionful (like the
known
$\kappa$-deformation) [14],
which makes them interesting for possible physical applications.
The other is that the standard q-deformation and the h-deformation
have been characterized as the only two deformations of $2\times 2$ matrices
with the central determinant [5, 6].

In this letter we propose to extend the h-deformation to the case of
supergroups
and discuss the simplest supergroup $GL(1|1)$.
We start with a superplane and its dual and follow the contraction method of
[8],
see also [15, 16].\\
Consider the $q$-deformed algebra of functions on the Manin superplane [13]
generated by $x^\prime , \theta^\prime$ with the relations
$x^\prime \theta^\prime = q ~\theta^\prime x^\prime$ and $(\theta^{\prime})^2 =
0$.

We introduce new coordinates  $x, \theta$ by
\begin{equation}
\pmatrix{x^\prime \cr \theta^\prime \cr} =g \pmatrix{x \cr \theta \cr} \ ,
\end{equation}
where
\begin{equation}
g=\pmatrix{1&0\cr
h/(q\! - \!1)&1\cr}\ .
\end{equation}
Here, it is worth pointing out that, differently to the usual situation,
the parameter $h$ can not be an ordinary number.
In fact, observing that $h$ describes
the mixing between the even coordinate $x$ and the odd one $\theta$,
it has to be a dual odd (grassman) number.
Thus we have $ h^2 = 0$ and $h\theta = - \theta h$.\\
Next, taking the limit $q \rightarrow 1$ we obtain the following
exchange relations, which define the h-superlane $A_h$
\begin{equation}
x \theta = \theta x + h x^2 , ~~~ {\theta}^2 = -h \theta x\ .
\label{ah}
\end{equation}
Similarly, we get the dual (exterior) $h$-superlane $\Lambda_h$ as generated by
$y, \xi$
with the relations
\begin{equation}
\xi y = - y \xi , ~~~ {\xi}^2 = 0\ .
\label{lh}
\end{equation}

We now define the corresponding $h$-deformation of the supergroup $GL(1|1)$
as a matrix quantum supergroup $GL_h(1|1)$ generatead by
$a, \beta , \gamma , d$ modulo the relations
\begin{eqnarray}
a\beta &=&{\beta}a,
{}~~~a\gamma = \gamma a + h (a^2 + \gamma \beta - d a ),
\nonumber\\
\beta d &=& d\beta   , ~~~\gamma d  =d \gamma + h (d^2 - \gamma \beta - d a ),
\nonumber\\
\beta ^2 &=& 0 , ~~~\gamma ^2 = h(d \gamma - \gamma a ),
\nonumber\\
\beta \gamma &=& - \gamma \beta -h(\beta a - d \beta ),
{}~~~a d  = d a + h(\beta a - d \beta).
\label{hcom}
\end{eqnarray}
These relations are obtained from the requirement that $A_h$ and $\Lambda_h$
have to be covariant under the (usual matrix) left coaction
\begin{equation}
\delta (x) = a \otimes x +\beta \otimes \theta ,
{}~~~ \delta(\theta )= \gamma\otimes x + d \otimes\theta
\end{equation}
and assuming that $\beta, \gamma$ anticommute with $\theta, y$ and $h$.

An interesting feature is that the first column of the matrix $T \in GL_h(1|1)$
is not isomorphic to $A_h$, though this is so in the case of q-deformations.

The relations (\ref{hcom}) can be alternatively obtained by performing the
similarity transformation introduced in [8]
\begin{equation}
T' = g T g^{-1}      \ ,
\end{equation}
which in our case reads
\begin{equation}
\pmatrix{a'&{\beta}'\cr
{\gamma}'&d'\cr}=
\pmatrix{a - (q\! - \!1)^{-1}\beta h&\beta\cr
(q\! - \!1)^{-1} h a +\gamma - (q\! - \!1)^{-1} d h&(q\! -\!1)^{-1} h \beta +
d\cr}  \ ,
\label{sim}
\end{equation}
and next the limit $q \rightarrow 1$.
Here, $a', {\beta}', {\gamma}', d'$
are generators of $GL_q(1|1)$, which satisfy the following commutation
relations:
\begin{eqnarray}
a'{\beta}' &=&q{\beta}'a',~~~~~~d'{\beta}'~= q{\beta}'d'~~~, \nonumber\\
a'{\gamma}' &=&q{\gamma}'a',~~~~~~d'{\gamma}'~= q{\gamma}'d'~~~, \nonumber\\
{{\beta}'}^2 &=&0,~{{\gamma}'}^2 =0,~{\beta}'{\gamma}'= -{\gamma}'{\beta}'~~,
\nonumber
\label{tcom}
\end{eqnarray}
\begin{equation}
ad-da=-(q-{1\over q}){\beta}'{\gamma}'\ .\nonumber
\label{qcom}
\end{equation}
Substituting (\ref{sim})
into (\ref{qcom}) we arrive at the set of relations (\ref{hcom}) above.

Hopf algebra structure is given by the usual
co-product and co-unit
\begin{equation}
{\bigtriangleup}(T_{ij}) = T_{ik} \otimes T_{kj} , ~~ i,j = 1,2,
\end{equation}
\begin{equation}
\varepsilon(T_{ij}) = \delta_{ij} ,
\end{equation}
Introducing the formal inverses $a^{-1}$ and $d^{-1}$,
and using the Gauss decomposition as in [11, 12] we find the antipode
\begin{equation}
S(T)=\pmatrix{a^{-1}+a^{-1}\beta d^{-1}\gamma a^{-1}&-a^{-1}\beta d^{-1}\cr
-d^{-1}\gamma a^{-1}&d^{-1}+d^{-1}\gamma a^{-1} \beta d^{-1}\cr}
\end{equation}
and the quantum super determinant (Berezinian) as
\begin{equation}
D_h = a (d-{\gamma}a^{-1}\beta)^{-1}\ .
\end{equation}
(Note, that these formulae are in fact independent of the relations
(\ref{tcom})). It is easy to see that $D_h$ can be rewritten equivalently as
\begin{eqnarray}
D_h&=&(d-{\gamma}a^{-1}\beta)^{-1} a \nonumber\\
&=&ad^{-1}-{\beta}d^{-1}{\gamma}d^{-1} \\
&=&d^{-1}a-d^{-1}{\beta}d^{-1}\gamma ~.\nonumber
 \end{eqnarray}
Moreover, it can be verified that $D$
has the `multiplicative property'
\begin{equation}
\Delta (D_h ) = D_h \otimes D_h \ .
\end{equation}
and that $D_h$ belongs to the centre of the algebra
\begin{equation}
TD=DT\ .
\end{equation}
Therefore by imposing the relation $D_h=1$, we may define
a deformation $SL_h(1|1)$.

The R-matrix for the supergroup $GL_h(1|1)$ can be
obtained from the R-matrix of $GL_q(1|1)$ by the transformation
\begin{equation}
R_h = (g \otimes g)^{-1} R_q (g \otimes g)
\end{equation}
Here [11]
\begin{equation}
R_q=\pmatrix{q&0&0&0\cr
0&1&0&0\cr
0&q - q^{-1}&1&0\cr
0&0&0&q^{-1}\cr}\ ,
\label{rq}
\end{equation}
and it is assumed that  $\otimes $ is graded  (we adhere to the conventions
of [12]) ie.
\begin{eqnarray}
(g_1)^{ab}_{cd} &=& (g\otimes I)^{ab}_{cd} =
(-1)^{c(b+d)} g^a_c {\delta}^b_d ~,\nonumber\\
(g_2)^{ab}_{cd} &=& (I\otimes g)^{ab}_{cd} =
(-1)^{a(b+d)} g^a_d {\delta}^a_c \ .
\end{eqnarray}
As a result, we obtain the following R-matrix $R_h$
\begin{equation}
R_h=\pmatrix{1&0&0&0\cr
-h&1&0&0\cr
h&0&1&0\cr
0&h&h&1\cr}\ ,
\label{rh}
\end{equation}
which satisfies the graded Yang-Baxter equation
\begin{equation}
R_{12}R_{13}R_{23}=R_{23}R_{13}R_{12}\ ,
\end{equation}
where
\begin{eqnarray}
(R_{12})^{i_1 i_2 i_3}_{j_1 j_2 j_3}&=&R^{i_1 i_2}_{j_1 j_2}
{\delta}^{i_3}_{j_3}~~,\nonumber\\
(R_{13})^{i_1 i_2 i_3}_{j_1 j_2 j_3}&=&(-1)^{i_2 (i_3 +j_3)}
R^{i_1 i_3}_{j_1 j_3} {\delta}^{i_2}_{j_2}~~,\label{gradr}\\
(R_{23})^{i_1 i_2 i_3}_{j_1 j_2 j_3}&=&(-1)^{i_1 (i_2 +i_3 +j_2 +j_3)}
R^{i_2 i_3}_{j_2 j_3} {\delta}^{i_1}_{j_1}\ , \nonumber
\end{eqnarray}
with $i_1, i_2, i_3, j_1, j_2, j_3 = 0, 1$.
Also,  $\hat R_h = {\cal P} R_h$, where
${\cal P}^{i_1 i_2}_{j_1 j_2}= (-1)^{i_1 i_2} {\delta}^{i_1}_{j_2}
{\delta}^{i_2}_{j_1}$
is the super permutation matrix, satisfies the graded braid equation
\begin{equation}
\hat R_{12}\hat R_{23}\hat R_{12}=\hat R_{23}\hat R_{12}\hat R_{23}\ ,
\end{equation}
with the grading again given by (\ref{gradr}).
Moreover, $\hat R_h$ satisfies
$$(\hat R_{h})^2 = I $$
and thus has two eigenvalues $\pm 1$.

We note that the commutation relations (\ref{hcom}) can be expressed by
\begin{equation}
R_h T_1 T_2 = T_2 T_1 R_h
\end{equation}
where $T_1$ and $T_2$ are graded in the same manner as $g_1$ and $g_2$
respectively.\\
Also, let us denote by $P_+$ and $P_-$ the projections onto the eigenspaces
$\pm 1$
of $\hat R_h$, and by $A$ and $\Lambda$ the quotients of algebras generated by
$x, \theta$ and $\xi , y$ modulo the ideal generated by ${\rm Ker}P_-$ and
${\rm Ker}P_+$, respectively.
Then $A$ and $\Lambda$ are isomorphic to $A_h$ defined by (\ref{ah})
and $\Lambda_h$ defined by (\ref{lh}), respectively.

We now pass to the differential calculus on the $h$-superplane $A_h$.
For this we insert the following solution
\begin{equation}
B = F = C = {\hat R}_h
\end{equation}
in the Wess-Zumino formulae $(2.1, 2.27, 2.28, 2.29$  and   $2.36)$ [10].
This yields:
\begin{eqnarray}
x \theta &=& \theta x + h x^2 , ~~~ {\theta}^2 = -h \theta x ,\nonumber \\
\xi y &=& - y \xi , ~~~ {\xi}^2 = 0 ,\nonumber \\
\partial_x \partial_{\theta} &=& \partial_{\theta} \partial_x , ~~~
\partial_{\theta}^2 = 0 ,\nonumber \\
x \xi &=& \xi x , ~~~\theta y = -y \theta -h \xi \theta -h y x ,\nonumber \\
x y &=& y x + h \xi x , ~~~ \theta \xi = \xi \theta -h \xi x ,\nonumber \\
\partial_x x &=& 1 + x \partial_x +h x \partial_{\theta} ,~~~
\partial_x \theta = \theta \partial_x -h x \partial_x
-h\theta\partial_{\theta},
\nonumber \\
\partial_{\theta} x &=& x \partial_{\theta} ,~~~ \partial_{\theta} \theta = 1 -
\theta \partial_{\theta} - h x \partial_{\theta},\nonumber \\
\partial_{\xi} \xi &=& \xi \partial_{\xi} + h \xi \partial_y ,~~~
\partial_{\xi} y = y \partial_{\xi} - h \xi \partial_{\xi} - h y \partial_y ,
\nonumber \\
\partial_y \xi &=& \xi \partial_y ,~~~ \partial_y y = - y \partial_y - h \xi
\partial_y \ .
\end{eqnarray}

It can be checked that this calculus is unique, satisfies all the consistency
conditions (see [10]) and is covariant under the coaction of $GL_h(1|1)$.

A few remarks are in order.
Note that $\hat R_h$ obeys only the graded braid equation and not the
ungraded one although  $\hat R_q = {\cal P} R_q$, with $R_q$ given by
(\ref{rq}),
is quite special as it satisfies both the ungraded and graded braid equations.
Moreover, $R_q$ and $P \hat R_q = diag(1, 1, 1, -1)R_q$,
where $P$ is the usual permutation, obey both the ungraded and graded
Yang-Baxter
equations. Surprisingly, due to the odd character of $h$ this is also the case
for $R_h$.

It is worth to mention that for the deformations of $GL(1|1)$
introduced in [12] with two parameters $p, q$,
if $p\neq q$ then $R$ given by (1) in [12] satisfies only
the graded Yang-Baxter equation and $diag(1, 1, 1, -1)R$ satisfies only
the (ungraded) Yang-Baxter equation although $\hat R = {\cal P} R$ obeys
both the ungraded and graded braid equations.

We close with a comment that it will be interesting to investigate the dual
deformation
$U_h gl(1|1)$ with odd parameter $h$, its representations and relation
with $U_q gl(1|1)$. Also the question of an isomorphism, upto Drinfeld twist
of the coproduct, with the Lie superalgebra $Ugl(1|1)$ should be studied.
The work on these issues is in progress.

\vskip.5cm
It is a pleasure to thank Prof. J. Lukierski for discussion on
$h$-deformations.
\newpage
{\bf References}
\newcounter{00001}
\begin{list}
{[~\arabic{00001}~]}{\usecounter{00001}
\labelwidth=1cm\labelsep=.5cm}

\item Demidov E.E, Manin Yu.I, Mukhin E.E and and Zhdanovich D.V, {\it Prog.
	Theor. Phys. Suppl. {\bf 102}}(1990)203.
\item Ewen H, Ogievetsky O and Wess J, {\it Lett. Math. Phys.} {\bf 22}(1991)
  297.
\item Zakrzewski S, {\it Lett. Math. Phys.} {\bf 22}(1991)287.
\item Ohn Ch, {\it Lett. Math. Phys.} {\bf 25}(1992)85.
\item Kupershmidt B.A, {\it J. Phys. A: Math. Gen.} {\bf 25}(1992)
   L1239.
\item Karimipour V, {\it Lett. Math. Phys.} {\bf 30}(1994)87.
\item Aghamohammadi A, {\it Mod. Phys. Lett. Math.} {\bf A8}(1993)2607.
\item Aghamohammadi A, Khorrami M and Shariati A, {\it J. Phys. A: Math. Gen.}
{\bf 28}(1995)L225 and references therein.
\item Reshetikhin N.Yu, Takhtajan L.A and Faddeev L.D, {\it Len. Math. J.
} {\bf 1}(1990)193.
\item Wess J and Zumino B, {\it Nucl. Phys. Proc. Suppl.} {\bf B18}(1990)302.
\item Schmidke W.B, Vokos S.P and Zumino B, {\it Z. Phys.} {\bf C48}(1990)249.
\item Dabrowski L and Wang Lu Yu, {\it Phys. Lett.} {\bf B266}(1991)51.
\item Manin Yu. I., {\it Quantum groups and non-commutative geometry}, {\bf
CRM},
Universite de Montreal,1988.
\item Lukierski J, Ruegg H and Tolstoy V. N., {ICTP preprint 1995}
\item Celeghini E, Giachetti R, Sorace E and Tarlini M, {\it J. Math. Phys.}
{\bf 32}(1991)1159.
\item Lukierski J, Nowicki A, Ruegg H and Tolstoy V. N., {\it Phys. Lett.}
{\bf B264}(1991)331.
\item Woronowicz S.L, {\it Rep. Math. Phys.} {\bf 30}(1991)259.

\end{list}
\end{document}